\documentclass[reqno]{amsart}
\usepackage{amssymb,setspace}
\usepackage{ifpdf}
\ifpdf
 \usepackage[hyperindex,pagebackref]{hyperref}%
\else
 \expandafter\ifx\csname dvipdfm\endcsname\relax
 \usepackage[hypertex,hyperindex,pagebackref]{hyperref}%
 \else
 \usepackage[dvipdfm,hyperindex,pagebackref]{hyperref}%
 \fi
\fi
\allowdisplaybreaks[4]
\numberwithin{equation}{section}
\theoremstyle{plain}
\newtheorem{thm}{Theorem}[section]
\newtheorem{cor}{Corollary}[section]
\theoremstyle{remark}
\newtheorem{rem}{Remark}[section]
\DeclareMathOperator{\td}{d}

\begin{document}

\title[Complete monotonicity of $q$-polygamma functions]
{Complete monotonicity of functions involving the $\boldsymbol{q}$-trigamma and $\boldsymbol{q}$-tetragamma functions}

\author[F. Qi]{Feng Qi}
\address{Institute of Mathematics, Henan Polytechnic University, Jiaozuo City, Henan Province, 454010, China}
\email{\href{mailto: F. Qi <qifeng618@gmail.com>}{qifeng618@gmail.com}, \href{mailto: F. Qi <qifeng618@hotmail.com>}{qifeng618@hotmail.com}, \href{mailto: F. Qi <qifeng618@qq.com>}{qifeng618@qq.com}}
\urladdr{\url{http://qifeng618.wordpress.com}}

\begin{abstract}
Let $\psi_q(x)$, $\psi_q'(x)$, and $\psi_q''(x)$ for $q>0$ stand respectively for the $q$\nobreakdash-digamma, $q$\nobreakdash-trigamma, and $q$\nobreakdash-tetragamma functions. In the paper, the author proves along two different approaches that the functions $[\psi'_q(x)]^2+\psi''_q(x)$ for $q>1$ and $[\psi_{q}'(x)-\ln q]^2 +\psi''_{q}(x)$ for $0<q<1$ are completely monotonic on $(0,\infty)$. Applying these results, the author derives monotonic properties of four functions involving the $q$-digamma function $\psi_q(x)$ and two double inequalities for bounding the $q$-digamma function $\psi_q(x)$.
\end{abstract}



\thanks{This paper was typeset using \AmS-\LaTeX}

\maketitle

\section{Introduction}

Recall from~\cite[Chapter~XIII]{mpf-1993} that a function $f$ is said to be completely monotonic on an interval $I$ if $f$ has derivatives of all orders on $I$ and
\begin{equation}\label{CM-dfn}
0\le(-1)^{n}f^{(n)}(x)<\infty
\end{equation}
for $x\in I$ and $n\ge0$.
For more information about properties and applications of completely monotonic functions, please refer to~\cite[Chapter~XIII]{mpf-1993} and closely related references therein.
\par
It is common knowledge that the logarithmic derivative of the classical Euler's gamma function $\Gamma(x)$, denoted by $\psi(x)=\frac{\Gamma'(x)}{\Gamma(x)}$, the derivatives $\psi'(x)$ and $\psi''(x)$ are respectively called the digamma, trigamma and tetragamma functions. As a whole, the derivatives $\psi^{(i)}(x)$ for $i\in\mathbb{N}$ are respectively called polygamma functions.
\par
Recall from~\cite[pp.~493\nobreakdash--496]{aar} and~\cite[Section~1.10]{Basic-hypergeometric-series-2nd} that the $q$\nobreakdash-gamma function, the $q$\nobreakdash-analogue of the gamma function $\Gamma(x)$, is defined for $x>0$ by
\begin{equation}\label{q-gamma-0-1}
\Gamma_q(x)=
\begin{cases}
(1-q)^{1-x}\prod\limits_{i=0}^\infty\dfrac{1-q^{i+1}}{1-q^{i+x}},& 0<q<1\\
(q-1)^{1-x}q^{\binom{x}2}\prod\limits_{i=0}^\infty\dfrac{1-q^{-(i+1)}}{1-q^{-(i+x)}},& q>1\\
\Gamma(x),& q=1
\end{cases}
\end{equation}
and satisfies
\begin{equation}\label{q-gamma-q-1/q}
\Gamma_q(x)=q^{\binom{x-1}2}\Gamma_{1/q}(x).
\end{equation}
The $q$-digamma function $\psi_q(x)$, the $q$\nobreakdash-analogue of the digamma function $\psi(x)$, is defined for $q>0$ and $x>0$ by
\begin{equation}
  \psi_q(x)=\begin{cases}\dfrac{\Gamma_q'(x)}{\Gamma_q(x)},&q\ne1,\\\psi(x),&q=1.
  \end{cases}
\end{equation}
The functions $\psi_q^{(k)}(x)$ for $k\in\mathbb{N}$, the $q$\nobreakdash-analogues of polygamma functions $\psi^{(k)}(x)$, are called $q$\nobreakdash-polygamma functions. In particular, the functions $\psi_q'(x)$ and $\psi_q''(x)$ are respectively called the $q$-trigamma and $q$-tetragamma functions.
From~\eqref{q-gamma-0-1}, we obtain that
\begin{enumerate}
\item
when $0<q<1$ and $x\in(0,\infty)$,
\begin{equation}\label{q-psi-dfn-q<1}
\begin{split}
\psi_q(x)&=-\ln(1-q)+(\ln q)\sum_{k=0}^\infty\frac{q^{k+x}}{1-q^{k+x}}\\
&=-\ln(1-q)+(\ln q)\sum_{k=1}^\infty\frac{q^{kx}}{1-q^{k}},
\end{split}
\end{equation}
\item
when $q>1$ and $x\in(0,\infty)$,
\begin{equation}\label{q-psi-dfn-q>1}
\psi_q(x)=-\ln(q-1)+(\ln q)\Biggl[x-\frac12-\sum_{i=1}^\infty\frac{q^{-ix}}{1-q^{-i}}\Biggr].
\end{equation}
\end{enumerate}
\par
The inequality
\begin{equation}\label{positivity}
\psi''(x)+[\psi'(x)]^2>0,\quad x>0
\end{equation}
was deduced or recovered in~\cite[Theorem~2.1]{batir-interest} and~\cite[Lemma~1.1]{batir-new}, employed in~\cite[Theorem~4.3]{alzer-grinshpan}, \cite[Theorem~2.1]{batir-interest}, and~\cite[Theorem~2.1]{batir-new},
and generalized in~\cite[Lemma~1.2 and Remark~1.3]{batir-jmaa-06-05-065}.
In~\cite{notes-best-simple-open-jkms.tex, notes-best-simple-cpaa.tex, AAM-Qi-09-PolyGamma.tex}, among other things, the inequality~\eqref{positivity} was generalized to complete monotonicity of the function $\psi''(x)+[\psi'(x)]^2$.
Recently, some other results in this field were published in~\cite{BAustMS-5984-RV.tex, SCM-2012-0142.tex} and closely related references therein.
\par
In~\cite[pp.~80--81, Lemma~4.6]{alzer-grinshpan}, the inequality~\eqref{positivity} was generalized as
\begin{equation}\label{positivity-q-psi}
  \bigl[\psi'_q(x)\bigr]^2+\psi''_q(x)>0
\end{equation}
for $q>1$ and $x>0$. Applying this inequality, the following conclusion was established in~\cite[p.~82, Theorem~4.8]{alzer-grinshpan}: for $q>1$ and $0<a<b\le\infty$, the double inequality
\begin{equation*}
  \exp\bigl\{\alpha\bigl[e^{\psi_q(x)}\bigl(\psi_q(x)-1\bigr)+1\bigr]\bigr\} \le\frac{\Gamma_q(x)}{\Gamma_q(x_0)} \le\exp\bigl\{\beta\bigl[e^{\psi_q(x)}\bigl(\psi_q(x)-1\bigr)+1\bigr]\bigr\}
\end{equation*}
is valid on $(a,b)$, with the best possible constant factors
\begin{equation*}
  \alpha=\begin{cases}
    Q_q(b),&\text{if $b<\infty$}\\0,&\text{if $b=\infty$}
  \end{cases}
\end{equation*}
and $\beta=Q_q(a)$, where $x_0=x_0(q)$ is the unique positive zero of the function $\psi_q$ for $0<q\ne1$, and
\begin{equation*}
  Q_q(x)=\begin{cases}
    \dfrac{\ln\Gamma_q(x)-\ln\Gamma_q(x_0)}{[\psi_q(x)-1]e^{\psi_q(x)}+1},&x\ne x_0,\\
    \dfrac1{\psi_q'(x_0)},&x=x_0.
  \end{cases}
\end{equation*}
For more information on the subject of $q$-gamma and $q$-polygamma functions, please refer to~\cite{aar, mon-element-exp-final.tex, Ismail-Muldoon-arXiv} and closely related references therein.
\par
The aim of this paper is to extend and generalize along two different approaches the inequality~\eqref{positivity-q-psi} to the following complete monotonicity.

\begin{thm}\label{psi-q-cm-thm}
The functions
\begin{equation}\label{psi-psi-2-functions-q}
  \mathcal{F}_q(x)=[\psi_{q}'(x)-\ln q]^2+\psi''_{q}(x), \quad 0<q<1
\end{equation}
and
\begin{equation}\label{psi-psi-2-functions}
  \mathfrak{F}_q(x)=\bigl[\psi'_q(x)\bigr]^2+\psi''_q(x), \quad q>1
\end{equation}
are completely monotonic on $(0,\infty)$.
\end{thm}

As consequences of Theorem~\ref{psi-q-cm-thm}, the following monotonicity and inequalities may be derived.

\begin{cor}\label{psi-q-ineq}
The functions
\begin{align}\label{phi-q-ge-1-dfn}
\phi_q(x)&=\psi_q(x)+\ln\biggl[\exp\dfrac{(\ln q)q^{x}}{q^{x}-1}-1\biggr],\quad q>1,\\
\label{phi-q-ge-1-dfn-1/p}
\varphi_q(x)&=\psi_{q}(x)+\ln\biggl(\exp\dfrac{\ln q}{q^{x}-1}-1\biggr),\quad q>1,\\
\label{phi-q<1-dfn}
  \Phi_q(x)&=\psi_{q}(x)-(\ln q)x+\ln\biggl[\exp\dfrac{(\ln q)q^{x}}{q^{x}-1}-1\biggr],\quad 0<q<1,\\
\label{phi-q<1-dfn-1/p}
\Theta_q(x)&=\psi_{q}(x)-(\ln q)x +\ln\biggl(\exp\dfrac{\ln q}{q^{x}-1}-1\biggr),\quad 0<q<1
\end{align}
are increasing on $(0,\infty)$. Consequently, the double inequalities
\begin{multline}\label{pai-q<1-upper}
\psi_{q}(1)+(\ln q)x-\ln\biggl[\exp\dfrac{(\ln q)q^{x}}{q^{x}-1}-1\biggr]<\psi_{q}(x)\\
<\ln\frac{\ln q}{q-1}+(\ln q)x-\ln\biggl[\exp\dfrac{(\ln q)q^{x}}{q^{x}-1}-1\biggr]
\end{multline}
for $0<q<1$ and
\begin{multline}\label{pai-q>1-ineq-double}
  \psi_{q}(1)-\ln q-\ln\biggl(\exp\dfrac{\ln q}{q^{x}-1}-1\biggr)<\psi_{q}(x)\\*
  <\ln\frac{\ln q}{q-1}-\frac{\ln q}2-\ln\biggl(\exp\dfrac{\ln q}{q^{x}-1}-1\biggr)
\end{multline}
for $q>1$ hold on $(0,\infty)$, where the constants $\psi_{q}(1)$ and $\ln\bigl(\frac{\ln q}{q-1}\bigr)$ in~\eqref{pai-q<1-upper} and the scalars $\psi_{q}(1)-\ln q$ and $\ln\frac{\ln q}{q-1}-\frac{\ln q}2$ in~\eqref{pai-q>1-ineq-double} are the best possible.
\end{cor}

\section{First proof of Theorem~\ref{psi-q-cm-thm}}

For $q>1$, using the expression~\eqref{q-psi-dfn-q>1} gives
\begin{equation}\label{psi-psi-2=deri}
  \frac{\bigl[\psi'_q(x)\bigr]^2+\psi''_q(x)}{(\ln q)^2} =1+\sum_{i=2}^\infty\sum_{j=1}^{i-1}c_j(q,x)c_{i-j}(q,x) -\sum_{i=1}^\infty(i-2)c_i(q,x),
\end{equation}
where
$$
c_i(q,x)=\frac{i\ln q}{1-q^{-i}}q^{-ix}>0.
$$
In order to prove complete monotonicity of the function~\eqref{psi-psi-2-functions}, by definition, it suffices to show
\begin{equation}\label{td-first-der}
(-1)^k\frac{\td{}^k}{\td x^k}\Biggl\{\frac{\bigl[\psi'_q(x)\bigr]^2+\psi''_q(x)}{(\ln q)^2}\Biggr\}\ge0,\quad k\ge0.
\end{equation}
For $k=0$, the inequality~\eqref{td-first-der} is equivalent to~\eqref{positivity-q-psi}. For $k\ge1$, it is equivalent to
\begin{equation}\label{partial-k-der}
(-1)^k\sum_{i=2}^\infty\sum_{j=1}^{i-1}\frac{\partial^k[c_j(q,x)c_{i-j}(q,x)]}{\partial x^k} \ge(-1)^k\sum_{i=1}^\infty(i-2)\frac{\partial^k[c_i(q,x)]}{\partial x^k}.
\end{equation}
In order to prove~\eqref{partial-k-der}, it is sufficient to verify
$$
(-1)^k\sum_{j=1}^{i-1}\frac{\partial^k[c_j(q,x)c_{i-j}(q,x)]}{\partial x^k} \ge(-1)^k(i-2)\frac{\partial^k[c_i(q,x)]}{\partial x^k},\quad i\ge3,
$$
that is,
\begin{multline*}
(-1)^k\sum _{j=1}^{i-1} \frac{(-1)^{k+1}i^{k} j (i-j) (\ln q)^{k+2} q^{i(1-x)+j}}{(q^j-1)(q^j-q^i)} \\ \ge(-1)^k(i-2)\frac{(-1)^{k+1}i^{k+1}(\ln q)^{k+1}q^{i(1-x)}}{1-q^i},\quad i\ge3,
\end{multline*}
which can be further simplified as
\begin{equation}\label{simplify-ineq}
\sum _{j=1}^{i-1} \frac{j (i-j)}{(q^j-1)(1-q^{i-j})} \le\frac{(i-2)i}{(\ln q)(1-q^i)},\quad i\ge3.
\end{equation}
Since
$$
\sum_{j=1}^{i-1}\frac{j(i-j)(q^{i-j}-q^j)}{(q^j-1)(q^{i-j}-1)}=0,
$$
we have
\begin{equation}\label{frac-j(i-j)(q^i-1)}
\sum_{j=1}^{i-1}\frac{j(i-j)(q^i-1)}{(q^j-1)(q^{i-j}-1)} =\sum_{j=1}^{i-1}\frac{j(i-j)(q^j+1)}{q^j-1}.
\end{equation}
The inequality
$$
e^\tau+1-\frac2\tau(e^\tau-1)=\sum_{k=2}^\infty\frac{k-1}{(k+1)!}\tau^k>0,\quad \tau=j\ln q>0
$$
implies that
$$
\frac{j(q^j+1)}{q^j-1}>\frac2{\ln q}.
$$
This leads to
\begin{equation*}
  \sum_{j=1}^{i-1}\frac{j(i-j)(q^j+1)}{q^j-1}>\frac2{\ln q}\sum_{j=1}^{i-1}(i-j)=\frac{i(i-1)}{\ln q} >\frac{i(i-2)}{\ln q}.
\end{equation*}
Substituting this inequality into~\eqref{frac-j(i-j)(q^i-1)} yields
\begin{equation}\label{frac-j(i-j)(q^i-1)op}
\sum_{j=1}^{i-1}\frac{j(i-j)(q^i-1)}{(q^j-1)(q^{i-j}-1)}>\frac{i(i-2)}{\ln q}
\end{equation}
which is equivalent to the strict inequality in~\eqref{simplify-ineq}. Thus, the inequality~\eqref{partial-k-der} is proved. Hence, the function $\mathfrak{F}_q(x)$ is completely monotonic on $(0,\infty)$.
\par
For $0<q<1$, taking the logarithm and differentiating on both sides of~\eqref{q-gamma-q-1/q} reveal
\begin{align}
\psi_q(x)&=(\ln q)\biggl(x-\frac32\biggr)+\psi_{1/q}(x),\label{q<1-2}\\
\psi_q'(x)&=\ln q+\psi_{1/q}'(x),\label{q<1-3}\\
\psi_q^{(k)}(x)&=\psi_{1/q}^{(k)}(x),\quad k\ge2.\label{q<1-4}
\end{align}
Therefore, when $0<q<1$, we have
\begin{equation}\label{Ffrak=cal-F(x)}
  \mathfrak{F}_{1/q}(x)=\bigl[\psi'_{1/q}(x)\bigr]^2+\psi''_{1/q}(x)
  =[\psi_{q}'(x)-\ln q]^2+\psi''_{q}(x)
  =\mathcal{F}_q(x)
\end{equation}
on $(0,\infty)$. As a result, the function $\mathcal{F}_q(x)$ is completely monotonic on $(0,\infty)$.
The first proof of Theorem~\ref{psi-q-cm-thm} is complete.

\section{Second proof of Theorem~\ref{psi-q-cm-thm}}

The formula~(1.11) in the paper~\cite{Ismail-Muldoon-119} and its corrected preprint~\cite{Ismail-Muldoon-arXiv} reads that
\begin{equation}\label{psi-q-int-eq}
\psi_q(x) = - \ln(1-q) -  \int_0^\infty\frac{e^{-xt}}{1-e^{-t}} d\gamma_q(t)
\end{equation}
for $0 < q < 1$ and $x >0$, where
\begin{equation}
  \gamma_q(t)=
\begin{cases}\displaystyle
-\ln q\sum_{k=1}^\infty\delta(t+k\ln q),& 0<q<1\\
t, & q=1
\end{cases}
\end{equation}
and $\delta(t)$ represents the Dirac delta function,
that is, $\td\gamma_q(t)$ is a discrete measure with positive masses $|\ln q|$ at the positive points $k|\ln q|$ for $k\in\mathbb{N}$.
Differentiating on both sides of the equation~\eqref{psi-q-int-eq} acquires
\begin{equation}\label{psi-q-der-eq}
\psi_q'(x)=\int_0^\infty\frac{te^{-xt}}{1-e^{-t}}\td\gamma_q(t), \quad 0<q<1.
\end{equation}
By the definition of $\gamma_q(t)$, we obtain that
\begin{equation}
\int_0^\infty e^{-xt}\td\gamma_q(t)=-\frac{q^x\ln q}{1-q^x}
\end{equation}
and
\begin{equation}
\int_0^\infty te^{-xt}\td\gamma_q(t)=\frac{q^x(\ln q)^2}{(1-q^x)^2}
\end{equation}
for $0<q<1$ and $x>0$.
It was presented in~\cite[p.~1245, Theorem~4.4, (4.15)]{mon-element-exp-final.tex}, \cite[Lemma~2.3 and Remark~2.1]{UPB-1635.tex}, \cite[Theorem~7.2, (7.5)]{bounds-two-gammas.tex}, and~\cite[p.~152, Theorem~4.22, (4.20)]{Wendel-Gautschi-type-ineq-Banach.tex} that, when $0<q<1$, the identity
\begin{equation}\label{n-s-ineq-g-gamma-equality}
\psi^{(k-1)}_q(x)-\psi^{(k-1)}_q(x+1)
=(\ln q)\frac{\td^{k-1}}{\td x^{k-1}}\biggl(\frac{q^{x}}{1-q^{x}}\biggr)
\end{equation}
is valid for $x\in(0,\infty)$ and $k\in\mathbb{N}$. Accordingly, it follows that
\begin{gather*}
\mathcal{F}_q(x)-\mathcal{F}_q(x+1)
=\bigl[\psi''_{q}(x)-\psi''_{q}(x+1)\bigr]\\
+\bigl[\psi_{q}'(x)-\psi_{q}'(x+1)\bigr] \bigl[\psi_{q}'(x)+\psi_{q}'(x+1)-2\ln q\bigr]\\
=(\ln q)\frac{\td^2}{\td x^2}\biggl(\frac{q^{x}}{1-q^{x}}\biggr)
+(\ln q)\frac{\td}{\td x}\biggl(\frac{q^{x}}{1-q^{x}}\biggr)
\biggl[2\psi_{q}'(x)-(\ln q)\frac{\td}{\td x}\biggl(\frac{q^{x}}{1-q^{x}}\biggr)-2\ln q\biggr]\\
=\frac{q^x (\ln q)^2}{(q^x-1)^2}
\biggl[2\psi_{q}'(x)-\frac{q^x (\ln q)^2}{(q^x-1)^2}-2\ln q\biggr]
-\frac{q^x(q^x+1) (\ln q)^3}{(q^x-1)^3}\\
=\frac{q^x (\ln q)^2}{(q^x-1)^2}
\biggl[2\psi_{q}'(x)-\frac{q^x (\ln q)^2}{(q^x-1)^2}-2\ln q
-\frac{(q^x+1) (\ln q)}{q^x-1}\biggr]\\
=\frac{2q^x (\ln q)^2}{(q^x-1)^2}
\biggl[\int_0^\infty\frac{te^{-xt}}{1-e^{-t}}\td\gamma_q(t)
-\frac12\int_0^\infty te^{-xt}\td\gamma_q(t)-\frac12\ln q
-\int_0^\infty e^{-xt}\td\gamma_q(t)\biggr]\\
=\frac{2q^x (\ln q)^2}{(q^x-1)^2}
\biggl[\int_0^\infty\biggl(\frac{t}{1-e^{-t}}-\frac12t-1\biggr)e^{-xt}\td\gamma_q(t) -\frac12\ln q\biggr]\\
=2\biggl[\int_0^\infty\biggl(\frac{t}{1-e^{-t}}-\frac12t-1\biggr)e^{-xt}\td\gamma_q(t) -\frac12\ln q\biggr]\int_0^\infty te^{-xt}\td\gamma_q(t).
\end{gather*}
By~\cite[p.~161, Theorem~12b]{widder} which reads that a function $f(x)$ is completely monotonic on $(0,\infty)$ if and only if it can be represented in the form
\begin{equation} \label{berstein-1}
f(x)=\int_0^\infty e^{-xt}\td\alpha(t),
\end{equation}
where $\alpha(t)$ is non-decreasing and the integral converges for $0<x<\infty$, and by the fact that the sum and product of finitely many completely monotonic functions are all still completely monotonic, we obtain that the difference $\mathcal{F}_q(x)-\mathcal{F}_q(x+1)$ is completely monotonic on $(0,\infty)$, that is,
\begin{align*}
0\le&(-1)^{k-1}[\mathcal{F}_q(x)-\mathcal{F}_q(x+1)]^{(k-1)}\\
&=(-1)^{k-1}\mathcal{F}_q^{(k-1)}(x)-(-1)^{k-1}\mathcal{F}_q^{(k-1)}(x+1)
\end{align*}
on $(0,\infty)$ for $k\in\mathbb{N}$. Consequently, by induction,
\begin{gather*}
(-1)^{k-1}\mathcal{F}_q^{(k-1)}(x)
\ge(-1)^{k-1}\mathcal{F}_q^{(k-1)}(x+1)
\ge(-1)^{k-1}\mathcal{F}_q^{(k-1)}(x+2)\\
\ge(-1)^{k-1}\mathcal{F}_q^{(k-1)}(x+3)
\ge\dotsm \ge(-1)^{k-1}\mathcal{F}_q^{(k-1)}(x+\ell)
\to0
\end{gather*}
as $\ell\to\infty$. This means that the function $\mathcal{F}_q(x)$ is completely monotonic on $(0,\infty)$.
\par
By~\eqref{Ffrak=cal-F(x)}, it is immediately deduced that the function $\mathfrak{F}_q(x)$ is completely monotonic on $(0,\infty)$.
The second proof of Theorem~\ref{psi-q-cm-thm} is complete.

\section{Proof of Corollary~\ref{psi-q-ineq}}

\subsection{}
Substituting~\eqref{q<1-3} and~\eqref{q<1-4} into~\eqref{n-s-ineq-g-gamma-equality}, we easily procure that
\begin{equation*}
\psi^{(k-1)}_{1/q}(x)-\psi^{(k-1)}_{1/q}(x+1)
=(\ln q)\frac{\td^{k-1}}{\td x^{k-1}}\biggl(\frac{q^{x}}{1-q^{x}}\biggr)
\end{equation*}
for $0<q<1$, $k\ge2$, and $x>0$. Replacing $q$ by $\frac1q$ in the above identity leads to
\begin{equation}\label{n-s-ineq-g-gamma-q>1}
\psi^{(k-1)}_{q}(x+1)-\psi^{(k-1)}_{q}(x)
=(\ln q)\frac{\td^{k-1}}{\td x^{k-1}}\biggl(\frac1{q^{x}-1}\biggr)
\end{equation}
for $q>1$ and $k\ge2$.
Similarly, combining~\eqref{q<1-2} with~\eqref{n-s-ineq-g-gamma-equality}, we obtain
\begin{equation}\label{n-s-ineq-g-gamma-equality-psi}
\psi_q(x+1)-\psi_q(x)=\frac{(\ln q)q^{x}}{q^{x}-1}
\end{equation}
for $q>1$ and $x>0$.
\par
When $0<q\ne1$, we have
\begin{align*}
\exp[\phi_q(x)]&=\exp[\psi_q(x)]\biggl[\exp\dfrac{(\ln q)q^{x}}{q^{x}-1}-1\biggr]\\
&=\exp\biggl[\psi_q(x)+\dfrac{(\ln q)q^{x}}{q^{x}-1}\biggr]-\exp[\psi_q(x)]\\
&=\exp[\psi_q(x+1)]-\exp[\psi_q(x)],
\end{align*}
in which the last equality is obtained by~\eqref{n-s-ineq-g-gamma-equality} and~\eqref{n-s-ineq-g-gamma-equality-psi}. Further differentiating gives
\begin{align*}
  \frac{\td\exp[\phi_q(x)]}{\td x}&=\psi_q'(x+1)\exp[\psi_q(x+1)]-\psi_q'(x)\exp[\psi_q(x)],\\
  \frac{\td\{\psi_q'(x)\exp[\psi_q(x)]\}}{\td x}&=\bigl\{\psi_q''(x)+[\psi_q'(x)]^2\bigr\}\exp[\psi_q(x)].
\end{align*}
The complete monotonicity of the function $\mathfrak{F}_q(x)$, obtained in Theorem~\ref{psi-q-cm-thm}, implies the positivity of $\mathfrak{F}_q(x)$ on $(0,\infty)$, that is, the inequality~\eqref{positivity-q-psi}. Hence, the function $\psi_q'(x)\exp[\psi_q(x)]$ is increasing, and so $\frac{\td\exp[\phi_q(x)]}{\td x}>0$ for $q>1$ on $(0,\infty)$. As a result, the function $\exp[\phi_q(x)]$, and thus $\phi_q(x)$, is increasing for $q>1$ on $(0,\infty)$.
\par
For $0<q\ne1$, we also have
\begin{align*}
  \exp[\Phi_q(x)]&=\exp[\psi_q(x)-(\ln q)x]\biggl[\exp\dfrac{(\ln q)q^{x}}{q^{x}-1}-1\biggr]\\
  &=\exp\biggl[\psi_q(x)+\dfrac{(\ln q)q^{x}}{q^{x}-1}-(\ln q)x\biggr]-\exp[\psi_q(x)-(\ln q)x]\\
  &=\exp[\psi_q(x+1)-(\ln q)x]-\exp[\psi_q(x)-(\ln q)x],
\end{align*}
in which the last equality is deduced by using~\eqref{n-s-ineq-g-gamma-equality} and~\eqref{n-s-ineq-g-gamma-equality-psi}. A direct computation gives
\begin{align*}
\frac{\td\exp[\Phi_q(x)]}{\td x}&=[\psi_q'(x+1)-\ln q]\exp[\psi_q(x+1)-(\ln q)x]\\
&\quad-[\psi_q'(x)-\ln q]\exp[\psi_q(x)-(\ln q)x]
\end{align*}
and
\begin{align*}
  &\quad\frac{\td\{[\psi_q'(x)-\ln q]\exp[\psi_q(x)-(\ln q)x]\}}{\td x}\\
  &=\bigl\{\psi_q''(x)+[\psi_q'(x)-\ln q]^2\bigr\}\exp[\psi_q(x)-(\ln q)x].
\end{align*}
The complete monotonicity of the function $\mathcal{F}_q(x)$, obtained in Theorem~\ref{psi-q-cm-thm}, implies the positivity of $\mathcal{F}_q(x)$ on $(0,\infty)$. Hence, the function
$$
\{[\psi_q'(x)-\ln q]\exp[\psi_q(x)-(\ln q)x]\}
$$
is increasing, and so $\frac{\td\exp[\Phi_q(x)]}{\td x}>0$ for $0<q<1$ on $(0,\infty)$. As a result, the function $\exp[\Phi_q(x)]$, and thus $\Phi_q(x)$, is increasing for $0<q<1$ on $(0,\infty)$.
\par
For $0<q<1$, utilizing~\eqref{q<1-2} yields
\begin{align*}
  \Phi_q(x)&=(\ln q)\biggl(x-\frac32\biggr)+\psi_{1/q}(x)-(\ln q)x+\ln\biggl[\exp\dfrac{(\ln q)q^{x}}{q^{x}-1}-1\biggr]\\
  &=\psi_{1/q}(x)+\ln\biggl[\exp\dfrac{(\ln q)q^{x}}{q^{x}-1}-1\biggr]-\frac32\ln q\\
  &=\psi_{1/q}(x)+\ln\biggl[\exp\dfrac{\ln (1/q)}{(1/q)^{x}-1}-1\biggr]+\frac32\ln\frac1q.
\end{align*}
Replacing $\frac1q$ by $q$ in the above equation leads to increasing monotonicity of the function $\varphi_q(x)$ for $q>1$ on $(0,\infty)$.
\par
For $q>1$, employing~\eqref{q<1-2} gives
\begin{align*}
\phi_q(x)&=\psi_{1/q}(x)-\biggl(\ln\frac1q\biggr)\biggl(x-\frac32\biggr) +\ln\biggl[\exp\dfrac{(\ln q)q^{x}}{q^{x}-1}-1\biggr]\\
&=\psi_{1/q}(x)-\biggl(\ln\frac1q\biggr)x +\ln\biggl[\exp\dfrac{\ln (1/q)}{(1/q)^{x}-1}-1\biggr] +\frac32\ln\frac1q.
\end{align*}
Substituting $\frac1q$ by $q$ in the above equation reduces to increasing monotonicity of the function $\Theta_q(x)$ for $0<q<1$ on $(0,\infty)$.

\subsection{}
From~\eqref{q-psi-dfn-q<1}, it is easy to see that
\begin{equation}\label{lim-x-to-infty-psi-q(x)}
  \lim_{x\to\infty}\psi_q(x)=-\ln(1-q),\quad q\in(0,1).
\end{equation}
In addition, for $q\in(0,1)$, we have
\begin{equation*}
  \lim_{x\to\infty}\biggl\{\ln\biggl[\exp\dfrac{(\ln q)q^{x}}{q^{x}-1}-1\biggr]-(\ln q)x\biggr\} =\lim_{x\to\infty}\ln\Biggl[\frac{\exp\frac{(\ln q)q^{x}}{q^{x}-1}-1}{q^x}\Biggr]=\ln\biggl(\ln\frac1q\biggr).
\end{equation*}
Combining them with increasing monotonicity of $\Phi_q(x)$, we acquire
$$
\psi_{q}(x)-(\ln q)x+\ln\biggl[\exp\dfrac{(\ln q)q^{x}}{q^{x}-1}-1\biggr] <\ln\biggl(\ln\frac1q\biggr)-\ln(1-q)=\ln\biggl(\frac{\ln q}{q-1}\biggr)
$$
for $q\in(0,1)$ on $(0,\infty)$, which is equivalent to the right-hand side of~\eqref{pai-q<1-upper}.
\par
From~\eqref{n-s-ineq-g-gamma-equality} for $k=1$ and the limit
\begin{equation}\label{lim-psi-0}
  \lim_{x\to0^+}\psi_{q}(x)=-\infty
\end{equation}
for $0<q<1$, it follows that
\begin{equation}\label{psi+q-lim}
\begin{aligned}
\psi_{q}(x)+\ln\biggl[\exp\dfrac{(\ln q)q^{x}}{q^{x}-1}-1\biggr] &=\ln\exp[\psi_{q}(x)]+\ln\biggl[\exp\dfrac{(\ln q)q^{x}}{q^{x}-1}-1\biggr]\\
&=\ln\biggl\{\biggl[\exp\dfrac{(\ln q)q^{x}}{q^{x}-1}-1\biggr]\exp[\psi_{q}(x)]\biggr\} \\
&=\ln\biggl\{\exp\biggl[\dfrac{(\ln q)q^{x}}{q^{x}-1}+\psi_{q}(x)\biggr]-\exp[\psi_{q}(x)]\biggr\} \\
&=\ln\bigl\{\exp\bigl[\psi_{q}(x+1)\bigr]-\exp[\psi_{q}(x)]\bigr\}
\to\psi_{q}(1)
\end{aligned}
\end{equation}
as $x\to0^+$, so the limit $\lim_{x\to0^+}\Phi_q(x)=\psi_{q}(1)$ for $0<q<1$. Combining this with increasing monotonicity of $\Phi_q(x)$, we acquire $\Phi_q(x)>\psi_{q}(1)$ for $0<q<1$ on $(0,\infty)$, that is,
$$
\psi_{q}(x)-(\ln q)x+\ln\biggl[\exp\dfrac{(\ln q)q^{x}}{q^{x}-1}-1\biggr]>\psi_{q}(1),\quad 0<q<1, \quad x\in(0,\infty).
$$
The required inequality in the left-hand side of~\eqref{pai-q<1-upper} is proved.
\par
For $0<q<1$, by using~\eqref{n-s-ineq-g-gamma-equality} for $k=1$, we have
\begin{equation}\label{Theta-q(x)-to0}
\begin{split}
\Theta_q(x)&=\psi_{q}(x)+\frac{(\ln q)q^{x}}{q^{x}-1}-\frac{(\ln q)q^{x}}{q^{x}-1}-(\ln q)x +\ln\biggl(\exp\dfrac{\ln q}{q^{x}-1}-1\biggr)\\
&=\psi_{q}(x+1)+\ln\frac{\exp[{\ln q}/{(q^{x}-1)}]-1}{q^x\exp[{(\ln q)q^{x}}/{(q^{x}-1})]}\\
&\to\psi_{q}(1)-\ln q
\end{split}
\end{equation}
as $x\to0^+$. Since $\Theta_q(x)$ is increasing on $(0,\infty)$, we obtain that
$$
\psi_{q}(x)-(\ln q)x +\ln\biggl(\exp\dfrac{\ln q}{q^{x}-1}-1\biggr)>\psi_{q}(1)-\ln q
$$
for $0<q<1$ on $(0,\infty)$, which is included in the left-hand side of~\eqref{pai-q<1-upper}.
\par
For $q>1$, by a similar argument as in~\eqref{Theta-q(x)-to0}, we have
\begin{equation*}
\varphi_q(x)=\psi_{q}(x+1)+\ln\frac{\exp[{\ln q}/{(q^{x}-1)}]-1}{\exp[{(\ln q)q^{x}}/{(q^{x}-1})]} \to\psi_{q}(1)-\ln q
\end{equation*}
as $x\to0^+$. Furthermore, from~\eqref{q<1-2} and~\eqref{lim-x-to-infty-psi-q(x)}, it follows that
\begin{align*}
\varphi_q(x)&=\psi_{1/q}(x)+(\ln q)\biggl(x-\frac32\biggr) +\ln\biggl(\exp\dfrac{\ln q}{q^{x}-1}-1\biggr)\\
&=\psi_{1/q}(x)-\frac{3\ln q}2+\ln\biggl[q^x\biggl(\exp\dfrac{\ln q}{q^{x}-1}-1\biggr)\biggr]\\
&\to-\ln\biggl(1-\frac1q\biggr)-\frac{3\ln q}2+\ln(\ln q)\\
&=\ln\frac{\ln q}{q-1}-\frac{\ln q}2
\end{align*}
as $x\to\infty$. Hence, we have
\begin{equation*}
\psi_{q}(1)-\ln q<\varphi_q(x)<\ln\frac{\ln q}{q-1}-\frac{\ln q}2
\end{equation*}
for $q>1$ on $(0,\infty)$. The required double inequality~\eqref{pai-q>1-ineq-double} is proved.
\par
For $q>1$, the limits~\eqref{lim-psi-0} and~\eqref{psi+q-lim} are still valid, then $\lim_{x\to0^+}\phi_q(x)=\psi_{q}(1)$ for $q>1$, and so $\phi_q(x)>\psi_{q}(1)$ which is contained in the left-hand side of~\eqref{pai-q>1-ineq-double}. The proof of Corollary~\ref{psi-q-ineq} is complete.

\section{Remarks}

Finally we remark some closely related concerns.

\begin{rem}
It is easy to see that
\begin{equation}\label{phi(x)-property-psi}
\lim_{q\to1}\phi_q(x)=\lim_{q\to1}\varphi_q(x)=\lim_{q\to1} \Phi_q(x)=\lim_{q\to1}\Theta_q(x)=\psi(x)+\ln\bigl(e^{1/x}-1\bigr).
\end{equation}
It is known~\cite{property-psi.tex-Hacet, property-psi-ii-Munich.tex} that the function~\eqref{phi(x)-property-psi} is increasing and concave on $(0,\infty)$ and that
\begin{equation}\label{batir-ineq-orig}
a-\ln\bigl(e^{1/x}-1\bigr)<\psi(x)<b-\ln\bigl(e^{1/x}-1\bigr)
\end{equation}
if and only if $a\le-\gamma$ and $b\ge0$, where $\gamma=0.577\dotsc$ stands for Euler-Mascheroni's constant. The inequality~\eqref{batir-ineq-orig} may also be derived from taking the limit $q\to1$ in~\eqref{pai-q<1-upper} or~\eqref{pai-q>1-ineq-double}. We conjecture that the functions $\phi_q(x)$, $\varphi_q(x)$, $\Phi_q(x)$, and $\Theta_q(x)$ should be concave on $(0,\infty)$.
\end{rem}

\begin{rem}
From the proof of Corollary~\ref{psi-q-ineq}, it may be obtained that the functions $\psi_q'(x)e^{\psi_q(x)}$ for $q\ge1$ and $[\psi_q'(x)-\ln q]e^{\psi_q(x)-(\ln q)x}$ for $0<q\le1$ are increasing on $(0,\infty)$.
In~\cite[Lemma~1.2]{batir-new} and~\cite[pp.~241--242]{egp}, it was obtained that the function $\psi'(x)e^{\psi(x)}$ is increasing and less than $1$ on $(0,\infty)$. Can one find bounds of the functions $\psi_q'(x)e^{\psi_q(x)}$ for $q\ge1$ and $[\psi_q'(x)-\ln q]e^{\psi_q(x)-(\ln q)x}$ for $0<q\le1$?
\end{rem}

\subsection*{Acknowledgements}
The author thanks anonymous referees for their comments on and helpful suggestions to the original version of this paper.

\end{document}